\documentclass[11pt]{article}
\usepackage{amsfonts,latexsym}
\usepackage{color}
\usepackage{mathptmx}
\usepackage{amsmath,amssymb}
\usepackage{amsthm}
\usepackage[mathscr]{euscript}
\usepackage{mathrsfs}

\newcommand{\I}{{\bf I}}

\newcommand{\rmd}{\mathrm{d}}

\newcommand{\sig}{\sigma}
\newcommand{\ph}{\varphi}

\newcommand{\Om}{\Omega}

\newcommand{\mc}[1]{\mathcal{#1}}
\newcommand{\mb}[1]{\mathbb{#1}}
\newcommand{\scr}[1]{\mathscr{#1}}

\newcommand{\bei}{\begin{itemize}}
\newcommand{\ei}{\end{itemize}}
\newcommand{\bg}{\begin}
\newcommand{\e}{\end}
\newcommand{\ed}{\end{document}}

\newcommand{\be}{\begin{enumerate}}
\newcommand{\ee}{\end{enumerate}}
\newcommand{\beq}{\begin{equation}}
\newcommand{\eeq}{\end{equation}}
\newcommand{\beqs}{\begin{equation*}}
\newcommand{\eeqs}{\end{equation*}}
\newcommand{\beqa}{\begin{eqnarray}}
\newcommand{\eeqa}{\end{eqnarray}}
\newcommand{\beqas}{\begin{eqnarray*}}
\newcommand{\eeqas}{\end{eqnarray*}}

\numberwithin{equation}{section}

\newcommand{\PrP}[2]{$\lk #1,\mb #2\rk$}

\newcommand{\lk}{\left(}
\newcommand{\rk}{\right)}

\newcommand{\OFP}{$\lk\Omega,\scr F,\mathbb{P}\rk$}
\newcommand{\OFPF}{$\lk\Omega,\scr F,\mathbb{P};\mathbb{F}\rk$}

\newtheorem{Th}{Theorem}[section]

\newtheorem{Prop}[Th]{Proposition}

\title{On the Minimal Entropy Martingale Measure\\ for L\'evy Processes}
\author{Andrii Andrusiv$^{\rm a}$ and Hans-J\"urgen Engelbert$^{\rm b}$\thanks{Corresponding author. Email: hans-juergen.engelbert@uni-jena.de}\\[2mm]
$^{\rm a}$Limited Liability Company Deloitte \& Touche, Business Center Prime\\
48, 50A Zhylyanska St., Kyiv, 01033,
Ukraine\\[2mm] 
$^{\rm b}$Institute of Mathematics, Friedrich SchillerUniversity of Jena\\
Ernst-Abbe-Platz 2, D-07743 Jena, Germany}

\date{}
\begin{document}
\thispagestyle{empty}
\maketitle
\begin{abstract}
\noindent In the present paper, a new and simple approach is provided for proving rigorously that for general L\'evy financial markets the minimal entropy martingale measure and the Esscher martingale measure coincide. The method consists in approximating the probability measure $\mb P$ by a sequence of L\'evy preserving probability measures $\mb P_n$ with exponential moments of all order. As a by-product, it turns out that the problem of finding the minimal entropy martingale measure for the L\'evy market is equivalent to the corresponding problem but for a certain one-step financial market. The existence of the Esscher martingale measure (and hence the minimal entropy martingale measure) will be characterized by using moment generating functions of the L\'evy process.  
\medskip  


\end{abstract}

\section{Introduction}\label{sec:intro}
We consider a geometric L\'evy market with asset price 
\beq\label{GLM}
S_t=S_0\exp(X_t),\quad t\in[0,T]\,,
\eeq
where the initial price $S_0>0$ is a constant, $T>0$ is a finite horizon, and \PrP X F\ is a L\'evy process on the filtered probability space \OFPF with characteristic triplet $(b,\sigma^2,\nu)$.
The interest rate is assumed to be equal to zero.

In general, the market is \textit{incomplete}: 
Only if $\lk X,\mb F\rk)$ is a Brownian motion \textit{or} a Poisson process, there exists a \textit{unique} equivalent martingale measure $\mb Q$. In all remaining cases, the set of equivalent martingale measures is either \textit{empty} or \textit{uncountable}. Therefore, the problem arises to choose suitable martingale measures $\mb Q$, if any, for pricing contingent claims. In some sense, $\mb Q$ should be ``close" to the physical probability measure $\mb P$.

Historically, two choices of equivalent martingale measures play an important role: The \textit{minimal entropy martingale measure} (abbreviated MEMM) and the \textit{Esscher martingale measure} (abbreviated EMM).

\smallskip
\noindent\textit{\textbf{Minimal Entropy Martingale Measure (MEMM)}} \ 
Let $\mathbb{Q}$ be a probability measure on $\lk\Omega,\scr F_T\rk$. For $t\in[0,T]$, the \textit{relative entropy} $I_t\lk\mathbb{Q},\mathbb{P}\rk$ of $\mathbb{Q}$ with respect to $\mathbb{P}$ on $\scr F_t$ is defined by 
\beq\label{RE}
\displaystyle I_t(\mathbb{Q},\mathbb{P}):=\begin{cases}\mathbb{E}_{\mathbb{P}}\left[\frac{\rmd\mathbb{Q}}{\rmd\mathbb{P}}\big|_{\scr F_t}\log\frac{\rmd\mathbb{Q}}{\rmd\mathbb{P}}\big|_{\scr F_t}\right],\quad &\textnormal{if } \ \mathbb{Q}|_{\scr F_t}\ll\mathbb{P}|_{\scr F_t},\\
\\+\infty,&\textnormal{otherwise}.
\end{cases}
\eeq
By $\scr M_a(S,\mb F)$ we denote the set of all \textit{absolutely continuous martingale measures} $\mb Q$ on $\scr F_T$: \  $\mathbb{Q}|_{\scr F_T}\ll\mathbb{P}|_{\scr F_T}\quad\mbox{and}\quad(S,\mb F)\quad\mbox{is a}\quad\mb Q\mbox{-martingale}$.
\bg{Def} $\mathbb{P}^*\in\scr M_a(S,\mb F)$ will be called the minimal entropy martingale measure (MEMM) if it satisfies
\beq\label{MEMM}
I_T(\mathbb{P}^*,\mathbb{P})=\inf_{\mathbb{Q}\in\scr M_a(S,\mb F)} I_T(\mathbb{Q} , \mathbb{P})\,.
\eeq
\e{Def}

The strong interest in the MEMM comes, in particular, from its close relation to \textit{portfolio optimization} in case of exponential utility (duality problem). See, e.g., Delbaen et al. \cite{Six-02} (the so-called six-author paper) or Kabanov \& Stricker \cite{KS02}.

Originally, relative entropy was introduced by Kullback--Leibler (sometimes also called \textit{Kullback--Leibler information number}).
Minimization problems for relative entropy with linear constraints have been studied in a pioneering paper by Csisz\'ar \cite{Cs75}. 

In his seminal paper, Frittelli \cite{Fr00} adapted Csisz\'ar's work to the needs of Mathematical Finance. In particular, in a general model Frittelli proved the existence of the minimal entropy martingale measure if the asset price processes are bounded, which is, however, a quite strong condition. He also proved that, if the minimal entropy martingale measure $\mb P^*$ exists, then $\mb P^*$ is equivalent to $\mb P$ on $\mc F_T$ if only there is some $\mb Q\in\scr M_a$ with finite entropy such that $\mb Q$ is equivalent to $\mb P$ on $\mc F_T$. 
For his general model, Frittelli has given an important characterization of the MEMM (cf. \cite[Theorem 2.3]{Fr00}).

\smallskip
\noindent\textit{\textbf{Esscher Martingale Measure (EMM)}} \ 
Esscher transformation and Esscher measure were introduced in Actuarial Mathematics by Fredrik Esscher \cite{E32} in 1932. Given a L\'evy process $(X,\mb F)$, for $\kappa\in\mathbb{R}$ with $\mathbb{E}\left[\exp\lk\kappa X_T\rk\right]<\infty$, let us define the probability density $Z_T^{\kappa}$,
\beqa\label{EM}
Z_T^{\kappa}:=\frac{\exp\lk\kappa X_T\rk}{\mathbb{E}\left[\exp\lk\kappa X_T\rk\right]}\,, 
\eeqa
and the probability measure $\mathbb{P}^{\kappa}$ with density $Z_T^{\kappa}$: \ $\rmd\mathbb{P}^{\kappa}=Z_T^{\kappa}\,\rmd \mathbb{P}$.
\bg{Def}\label{DEMM}
{\rm (i)} The probability measure $\mathbb{P}^{\kappa}$ is called Esscher measure.

{\rm (ii)} \ $\mathbb{P}^{\kappa}$ is called Esscher martingale measure for the geometric L\'evy market $(S,\mb F)$ (resp., for the linear L\'evy market $(X,\mb F)$) if the process $\lk S,\mb F\rk$ (resp., $(X,\mb F)$) is a martingale with respect to $\mathbb{P}^{\kappa}$.

{\rm (iii)} \ The EMM, if it exists, will be denoted by $\mb P^{E,g}$ (resp., $\mb P^E$).
\e{Def}
The two kinds of EMM for the geometric and linear L\'evy markets should be carefully distinguished. We shall return to this point later on. 

The EMM $\mathbb{P}^{E,g}$ (resp., $\mathbb{P}^E$) is obviously unique. However, the existence of an Esscher measure $\mathbb{P}^{\kappa}$ different from $\mb P$ requires finiteness of the exponential moment $\mathbb{E}\left[\exp\lk\kappa X_t\rk\right]$ for $t=T$ and hence for all $t\in[0,T]$. If there \textit{does} exist some finite exponential moments $\mathbb{E}\left[\exp\lk\kappa X_t\rk\right]$, then Esscher measures exist, but there \textit{need not} exist the EMM. 

Esscher measures have many useful properties. In particular, the Esscher transformation preserves the L\'evy property: \ $\lk X,\mb F\rk$ is again a $\mathbb{P}^{\kappa}$-L\'evy process.

Gerber \& Shiu \cite{GS94,GS96} first suggested to use the EMM for \textit{option pricing}. This idea was significantly developed, in particular, because of the duality between the problem of portfolio optimization and the minimization of relative entropy.

\smallskip
\noindent\textit{\textbf{Comparing EMM and MEMM}} \ 
In this subsection, we are going to discuss what is known from the literature about the relation between EMM and MEMM.
 
Under the condition of finite exponential moments in a neighbourhood of zero, Chan \cite{C99} gives arguments of \textit{heuristical nature} that the minimal entropy martingale measure is just the same as the EMM for a linear L\'evy market. This is realized by approaching the minimization problem \eqref{MEMM} restricting to \textit{deterministic} Girsanov parameters for the densities of $\mb Q$ with respect to $\mb P$ on $\mc F_T$. However, this would require a rigorous proof.

Fujiwara \& Miyahara \cite{FM03} give analytical conditions in terms of the characteristic triplet $(b,\sig,\nu)$  of the L\'evy process that the EMM exists and prove that the EMM is the MEMM. The conditions are: \ $\exists\,\kappa^*\in\mb R$ such that
\beqas
(1)&& \quad \int_\mb R \big|x\,\exp\lk\kappa^*\,x\rk-h(x)\big|\,\nu(\rmd x)<+\infty\,;\\
(2)&& \quad b+\sig\,\kappa^*+\int_\mb R \lk x\,\exp\lk\kappa^*\,x\rk-h(x)\rk\,\nu(\rmd x)=0\,.
\eeqas
In probabilistic terms this means:
\beqas
\hspace{-2,75cm}(1^\prime)&& \quad\mb E\left[\big|L_T\big|\exp\lk\kappa^*\,L_T\rk\right]<+\infty\,;\\
\hspace{-2,75cm}(2^\prime)&&\quad\mb E\big[L_T\,\exp(\kappa^*\,L_T)\big]=0\,,
\eeqas
which is, however, just the definition of the EMM for $(L,\mb F)$. Here $(L,\mb F)$ denotes the linear L\'evy process associated with the geometric L\'evy process $(S,\mb F)$, see Section \ref{Sec:Prel} for details. 

Esche \& Schweizer \cite{ES05} suggest a rigorous proof of the result given by Chan \cite{C99}. For this, they justify the minimization over deterministic Girsanov parameters (see Theorem A in \cite{ES05}). However, the solution of the infinite dimensional minimization problem for the MEMM, leading to the EMM, is demanding and seems a bit formal: The authors admit a ``heuristically derived recipe" and that the derivation of the candidate for the MEMM ``by partly formal arguments" (see \cite[Section 4]{ES05}).  Moreover, the proof of Theorem A is quite sophisticated, using the whole machinery of stochastic analysis for averaging arbitrary Girsanov parameters to get deterministic ones. Theorem B in \cite{ES05} extends the paper of Fujiwara \& Miyahara \cite{FM03}, in particular, to the multidimensional case. 

Hubalek \& Sgarra \cite{HS06} offer a proof that the MEMM is the EMM for the linear L\'evy market $(L,\mb F)$ associated with the geometric L\'evy market $(S,\mb F)$ (see Section \ref{Sec:Prel}). However, they rely on the main results of Esche and Schweizer \cite[Theorems A and B]{ES05}, so their proof cannot be considered as independent of \cite{ES05}, and there remain some unclear points in the proof. The novel idea is to replace the infinite-dimensional minimization problem by a one-dimensional one by suitably disturbing the density of the L\'evy measure for the minimal entropy martingale measure. 

In her thesis, Cawston \cite{Ca-Th10} states that MEMM and EMM coincide, as a consequence of her Proposition 2.7, referring to Fujiwara \& Miyahara \cite{FM03} for one part and to Hubalek \& Sgarra \cite{HS06} for the other.

Cawston \& Vostrikova \cite{C-V-3,C-V-1,C-V-2} discuss some properties of MEMM and Esscher measures for L\'evy models, as well as the preservation of the L\'evy property under f-minimal measures.

\smallskip 
The main objective of the present paper is to give a new and simple approach for a rigorous proof that the MEMM \ $\mb P^*$ and the EMM \ $\mb P^E$ for a linear L\'evy market in fact coincide. More precisely, if one of these probability measures exists, then there exists the other, and they are equal: \ $\mb P^*=\mb P^E$. The method consists in approximating the probability measure $\mb P$ by a sequence of L\'evy preserving probability measures $\mb P_n$ with exponential moments of all order. As a useful by-product, from this it can be derived that the problem of finding the minimal entropy martingale measure for the L\'evy market is equivalent to the corresponding problem but for a one-step financial market. As a particular result, denoting by $\mb P^*_n$ the MEMM with respect to $\mb P_n$, the entropy $I_T(\mathbb{P}^*_n,\mathbb{P})$ approximates the lower bound of entropy $I_T(\mathbb{Q},\mathbb{P})$ over different classes of absolutely continuous probability measures $\mb Q$, even in the case that the MEMM $\mb P^*$ with respect to $\mb P$ does not exist.

For simplification of the presentation, we shall focus on the case of one-dimensional L\'evy processes. However, we emphasize that our approach and the results stated in the present paper can be extended to arbitrary multidimensional L\'evy processes by working along similar lines.

In Section \ref{Sec:Prel}, we start with some definitions and notations. We continue discussing geometric L\'evy processes versus linear L\'evy processes. Then we recall the existence of the EMM provided that exponential moments of arbitrary order are finite. We proceed with a discussion of the no-arbitrage condition. Finally, we recall that the EMM, if it exists, is always the MEMM and give a simple proof of this fact.

In Section \ref{sec:App}, we describe the approximation procedure for the probability $\mb P$ specifying the L\'evy market by suitable L\'evy preserving probabilities $\mb P_n$. In Propositions \ref{identity1} -- \ref{identity3}, we study the limit behaviour of different kinds  of entropy taken with respect to $\mb P_n$.

In Section \ref{sec:con}, we state the main results of the present paper, among them Theorem \ref{identity4} and Theorem \ref{mainth} and their corollaries, which are derived on the basis of the results prepared in the foregoing sections.

In the Appendix, for the convenience of the reader, we collect some properties of moment generating functions.

\section{Preliminaries}\label{Sec:Prel}
\noindent\textit{\textbf{Some Definitions and Notations}} \ Let $(X,\mb F)$ be an arbitrary L\'evy process on \OFP.
\bg{Def}\label{Lpreserving}
A probability measure $\mb Q$ on $(\Om,\mc F)$ is called L\'evy preserving for $(X,\mb F)$ if $(X,\mb F)$ is a L\'evy process with respect to $\mb Q$.
\e{Def}
\bg{Def}\label{moment-condition} {\rm (i)} A probability measure $\mb Q$ on $(\Om,\mc F)$ is said to satisfy the moment condition for $X$ (denoted by $\mb Q\in{\rm MC}(X)$) if $\mb E_\mb Q\left[X_T\right]=0$.

{\rm (ii)} A probability measure $\mb Q$ on $(\Om,\mc F)$ is said to satisfy the local moment condition for $(X,\mb F)$ (denoted by $\mb Q\in{\rm MC^{loc}}(X,\mb F)$) if there exists a sequence $(\tau_k)$ of $\mb F$-stopping times such that $\{\tau_k=T\}\uparrow_{k\to\infty}\Om$ and $\mb E_\mb Q\left[X_{\tau_k\wedge T}\right]=0$ for all $k\ge1$.
\e{Def}
The following classes of probability measures will be considered:
\beqa
{\scr M}_a(X,\mb F)&=&\{\mb Q: \ \mb Q|_{\scr F_T}\ll\mb P|_{\scr F_T}, \ (X,\mb F) \mbox{ is a } \mb Q\mbox{-martingale}\}\,;\label{abs-mart}\\
{\scr M}_e(X,\mb F)&=&\{\mb Q: \ \mb Q|_{\scr F_T}\sim\mb P|_{\scr F_T}, \ (X,\mb F) \mbox{ is a } \mb Q\mbox{-martingale}\}\,;\label{equ-mart}\\
{\scr M}_f(X,\mb F)&=&\{\mb Q\in{\scr M}_a(X,\mb F): \ I_T(\mb Q,\mb P)<+\infty\}\,;\label{fin-ent}\\
{\scr M}_f(X,\mb F)&=&\{\mb Q\in{\scr M}_a(X,\mb F): \ \mb Q \mbox{ is L\'evy preserving}\}\,;\label{Levy-pres-mart}\\
{\scr M}_{efl}(X,\mb F)&=&{\scr M}_e(X,\mb F)\cap{\scr M}_f(X,\mb F)\cap{\scr M}_l(X,\mb F)\,;\label{efl-mart}\\
{\scr M}^{\rm loc}_a(X,\mb F)&=&\{\mb Q: \ \mb Q|_{\scr F_T}\ll\mb P|_{\scr F_T}, \ (X,\mb F) \mbox{ is a } \mb Q\mbox{-local martingale}\}\,;\label{abs-locmart}\\
\tilde{{\scr M}}_a(X,\mb F)&=&\{\mb Q: \ \mb Q|_{\scr F_T}\ll\mb P|_{\scr F_T}, \ \mb Q\in{\rm MC}(X)\}\,;\label{abs-mom}\\
\tilde{{\scr M}}^{\rm loc}_a(X,\mb F)&=&\{\mb Q: \ \mb Q|_{\scr F_T}\ll\mb P|_{\scr F_T}, \mb Q\in{\rm MC^{loc}}(X,\mb F)\}\,;\label{abs-loc-mom}\\
\tilde{{\scr M}}^{\rm loc,0}_a(X,\mb F)&=&\{\mb Q\in\tilde{{\scr M}}^{\rm loc}_a(X,\mb F): \ \mb E_\mb Q\left[|X_T|\right]<+\infty\}\,;\label{abs-loc-mom-int}\\
\hat{{\scr M}}^{\rm loc}_a(X,\mb F)&=&{\scr M}^{\rm loc}_a(X,\mb F)\cup\tilde{{\scr M}}^{\rm loc,0}_a(X,\mb F)\,.\label{abs-loc-hat}
\eeqa
The reader may wonder at this place why we have introduced classes of absolutely continuous probability measures as \eqref{abs-mom} -- \eqref{abs-loc-hat} satisfying some (local) moment condition instead of (local) martingale property. Although these classes are major extensions of ${\scr M}_a(X,\mb F)$ and ${\scr M}^{\rm loc}_a(X,\mb F)$, intuitively speaking, the minimal entropy even taken over the biggest class $\tilde{{\scr M}}^{\rm loc}_a(X,\mb F)$ should not be reduced. As a by-product of our approach, we shall provide a rigorous proof of this basic fact. The heuristical idea comes from Lemma \ref{LocMoment} below.   
\bg{Le}\label{LocMoment}
{\rm (i)} Suppose that $\mb Q\in{\scr M}^{\rm loc}_a(X,\mb F)$ is L\'evy preserving. Then $\mb Q\in{\scr M}_a(X,\mb F)$.

\smallskip
{\rm (ii)} Now suppose that $\mb Q\in\tilde{\scr M}_a^{\rm loc,0}(X,\mb F)$ is L\'evy preserving. Then $\mb Q\in{\scr M}_a(X,\mb F)$. 
\e{Le}
\proof (i) By assumption, $(X,\mb F)$ is a L\'evy process and a local martingale with respect to $\mb Q$. According to He et al. \cite[Theorem 11.46]{HWY92}, the process $(X,\mb F)$ is martingale with respect to $\mb Q$, hence $\mb Q\in{\scr M}_a(X,\mb F)$.

(ii) The condition $\mb Q\in\tilde{\scr M}_a^{\rm loc,0}(X,\mb F)$ yields that $\mb E_\mb Q\left[|L_T|\right]<+\infty$ and that there exists a sequence $(\tau_k)$ of $\mb F$-stopping times such that $\{\tau_k=T\}\uparrow_{k\to\infty}\Om$ and $\mb E_\mb Q\left[X_{\tau_k\wedge T}\right]=0$ for all $k\ge1$. It is clear that $\mb E_\mb Q\left[X_{t}\right]=t\,\mb E_\mb Q\left[X_{1}\right]$, for all $t\in[0,T]$. Defining $M_t:=X_t-t\,\mb E_\mb Q\left[X_1\right]$, $t\in[0,T]$, implies that $(M,\mb F)$ is a martingale starting from zero. Now it follows that the $\mb Q$-expectation of $X_{\tau_k\wedge T}-M_{\tau_k\wedge T}=(\tau_k\wedge T)\,\mb E_\mb Q\left[X_1\right]$ is zero, which yields $\mb E_\mb Q\left[(\tau_k\wedge T)\right]\,\mb E_\mb Q\left[X_1\right]=0$ for all $k\ge0$, hence $\mb E_\mb Q\left[X_1\right]=0$. Consequently, $\mb E_\mb Q\left[X_t\right]=0$ for all $t\in[0,T]$ and, as a L\'evy process, $(X,\mb F)$ is a $\mb Q$-martingal.\hfill$\Box$

\smallskip 
\noindent\textit{\textbf{Geometric L\'evy Processes versus Linear L\'evy Processes}} \ In the Introduction, as common in Mathematical Finance, we started from the geometric L\'evy Market \eqref{GLM},
$$
S_t=S_0\exp(X_t),\quad t\in[0,T],
$$
where $S_0>0$ and $(X,\mb F)$ is a L\'evy process. However, the asset price $(S,\mb F)$, being a strictly positive semimartingale, can be represented as \textit{stochastic exponential} of another semimartingale $(L,\mb F)$: \ $S_t=S_0\,\scr E(L)_t, \ t\in[0,T]$,  
where $L:=\scr L(S)$ is the \textit{stochastic logarithm} of $S$ (see Jacod \& Shiryaev \cite[Ch. II.8]{JS00}). Note that $(L,\mb F)$ is again a L\'evy process, with characteristic triplet $(\tilde{b},\tilde{\sig}^2,\tilde{\nu})$ given by 
\beqas
\tilde{b}&=&b+\sig^2/2+\int_\mb R \left[\left(e^x-1\right)\I_{\{|e^x-1|\le1\}}-x\,I_{\{|x|\le1\}}\right]\,\nu(\rmd x)\,,\\ 
\tilde{\sig}^2&=&\sig^2/2,\quad \tilde{\nu}(B)=\int_\mb R \I_B(e^x-1)\,\nu(\rmd x), \ B\in\scr{B}\lk \mb R\rk
\eeqas
(cf. \cite[Corollary II.8.16]{JS00}). It should be emphasized that the L\'evy measure $\tilde\nu$ of $(L,\mb F)$ is supported by $(-1,+\infty)$ (in the sense of $\tilde{\nu}((-\infty,-1])=0$) which is equivalent to the property that $L$ only admits jumps bigger than $-1$. Moreover, we have the following proposition. 
\bg{Prop}\label{lmm}
Let $\mb Q$ be a probability measure on $(\Om,\mc F_T)$.

{\rm (i)} \ $\mb Q$ is a local martingale measure for $(S,\mb F)$ if and only if $\mb Q$ is a local martingale measure for $(L,\mb F)$.

{\rm (ii)} \ Let $\mb Q$ be L\'evy preserving for $(L,\mb F)$. Then $\mb Q$ is a martingale measure for $(S,\mb F)$ if and only if $\mb Q$ is a martingale measure for $(L,\mb F)$. 
\e{Prop}
\proof The proof of (i) is straightforward and therefore omitted. For proving (ii), let $\mb Q$ be a martingale measure for $(S,\mb F)$ which is L\'evy preserving for $(L,\mb F)$. Using (i), we obtain that $(L,\mb F)$ is a local martingale with respect to $\mb Q$. Hence, by He et al. \cite[Theorem 11.46]{HWY92}, the $\mb Q$-local martingale $(L,\mb F)$, being a $\mb Q$-L\'evy process, must be a martingale with respect to $\mb Q$. Conversely, let $\mb Q$ be a L\'evy preserving martingale measure for $(L,\mb F)$. Since $S=S_0\,\scr{E}\lk L\rk$, where $\scr{E}\lk L\rk$ denotes the stochastic exponential of $L$, from Cont \& Tankov \cite[Proposition 8.23]{CT03} we obtain that $(S,\mb F)$ is a martingale with respect to $\mb Q$. \hfill $\Box$

\smallskip
Later we shall prove that a sufficient subclass of martingale measures $\mb Q$ for solving the minimization problem $\inf_{\mb Q\in\scr M_a(L,\mb F)} I_T(\mb Q,\mb P)$ is the class of all L\'evy preserving martingale measures $\mb Q$. As a result of Proposition \ref{lmm}, it is therefore equivalent to deal with either the geometric L\'evy market $(S,\mb F)$ or with the linear L\'evy market $(L,\mb F)$.  

In the following, we prefer to deal with a \textit{general} linear L\'evy market $(L,\mb F)$, which is much more convenient and even more general. For the sake of simplicity, the characteristic triplet of $(L,\mb F)$ is again denoted by $(b,\sig^2,\nu)$, with $b\in\mb R$, $\sig^2\ge0$, and now general L\'evy measure $\nu$. In what follows, the classes of probability measures \eqref{abs-mart} -- \eqref{abs-loc-hat} will only be used for $(X,\mb F)=(L,\mb F)$ and, for the sake of simplification of notation, the suffix ``$(L,\mb F)$" will be omitted. Thus, from now on, e.g., ${\scr M}_a$ always stands for ${\scr M}_a(L,\mb F)$.

\smallskip
\noindent\textit{\textbf{The EMM: \ Preliminary Step}} \ Let $(L,\mb F)$ be a L\'evy process with characteristic triplet $(b,\sig^2,\nu)$. In this subsection, we assume that $L_T$ has finite exponential moments of any order, i.e.,
$
I:=\{\kappa\in\mb R: \ \mathbb{E}[\exp(\kappa L_T)]<+\infty\}=\mb R
$.
We have introduced the Esscher measures $\mb P^\kappa$ by its density \eqref{EM}: \ $Z_T^{\kappa}:=\exp(\kappa L_T)/\mathbb{E}[\exp(\kappa L_T)]$.
Note that $\mb P^{\kappa_0}\in\scr M_e$ (the set of equivalent martingale measures) if and only if 
\beq\label{EMM}
\mb E_{\mb P^{\kappa_0}}\left[L_T\right]=\mb E_{\mb P}\left[L_T\,Z_T^{\kappa_0}\right]=0\,.
\eeq
In this case, $\mb P^{\kappa_0}$ is called the Esscher martingale measure. We now introduce the functions $\varphi_T$ and $\psi_T$ by
\beq\label{CF}
\varphi_T(\kappa)=\mb E_{\mb P}\big[\exp(\kappa L_T)],\quad \psi_T(\kappa)=\mb E_{\mb P}\big[L_T\exp(\kappa L_T)], \quad \kappa\in\mb R\,.
\eeq
It can easily be calculated that the Esscher martingale measure $\mb P^{\kappa_0}$ has the finite entropy 
\beq\label{EM-ent}
I_T(\mb P^{\kappa_0},\mb P)=-\log\varphi_T(\kappa_0)\,.
\eeq
\bg{Prop}\label{Ex-EMM}
Suppose that $L$ is not monotone and that $\varphi_T(\kappa)<\infty$ for all $\kappa\in\mb R$. Then there exists the EMM $\mb P^E$. 
\e{Prop}
\proof Under the condition of the proposition that the L\'evy process $(L,\mb F)$ is not monotone, it follows that $\mb P(\{L_T<0\})>0$ and $\mb P(\{L_T>0\})>0$. Hence, Proposition \ref{prop:psi} can be applied for $\xi=L_T$. Therefore, $\psi_T$ is a continuous and increasing real-valued function satisfying 
$$
\lim_{\kappa\rightarrow\infty}\psi_T(\kappa)=+\infty,\quad \lim_{\kappa\rightarrow-\infty}\psi_T(\kappa)=-\infty\,.
$$
Consequently, there exists $\kappa_0\in\mb R$ such that
$$
\mb E_{\mb P}\big[L_T\exp(\kappa_0 L_T)]=\psi_T(\kappa_0)=0\,,
$$
meaning that $\mb P^{\kappa_0}$ is the EMM $\mb P^E$. \hfill$\Box$

\smallskip
\noindent\textit{\textbf{No-Arbitrage Condition}} \ 
Looking for the minimal entropy martingale measure for $(L,\mb F)$, as a necessary hypothesis, requires that the set $\scr M_a$ is \textit{nonempty}, meaning that the market should have no arbitrage. We will now give a short and simple proof of the following lemma which was stated in Cherny \& Shiryaev \cite{CS02}. 
\bg{Le}\label{No-Arb}
The market is free of arbitrage if and only if $L$ is not monotonic, that is, $(L,\mb F)$ and $(-L,\mb F)$ are no subordinators.
\e{Le}
\proof Only the sufficiency of the condition is to show. Here is a brief outline, cf. the beginning of Section \ref{sec:App} and Proposition \ref{P_n:Levy-preserving} in case of $n=1$ for more details. Define $Y(x)=\exp(-\rho(x))$ where $\rho(x)=x^2$ if $|x|>1$ and equal to $0$ otherwise, and for $t\in[0,T]$ put
\beq\label{density}
Z_t=\scr E\big((Y-1)*(\mu-\mu^p)\big)_t=\exp\lk-\sum_{0<u\le t}\rho(\Delta L_u)+t\,\int_\mb R \lk 1-Y(x)\rk\,\nu(\rmd x)\rk\,,
\eeq
where $(Y-1)*(\mu-\mu^p)$ is the stochastic integral with respect to the compensated Poisson random measure $\mu-\mu^p$ of $L$ and $\scr E(S)$ denotes the stochastic exponential of a semimartingale $S$. Then $(Z,\mb F)$ is a bounded nonnegative martingale with expectation $1$. Define the probabilty measure $\mb P_1$ by $\rmd \mb P_1=Z_T\,\rmd \mb P$ which is obviously equivalent to $\mb P$. The probability measure $\mb P_1$ is L\'evy preserving (see Definition \ref{Lpreserving}) and the L\'evy measure $\nu_1$ of $(L,\mb F)$ with respect to $\mb P_1$ is $\rmd \nu_1=Y(x)\,\rmd\nu(x)$. It follows that $L_T$ has finite exponential moments of arbitrary order \textit{with respect to} $\mb P_1$ (cf. Sato \cite[Theorem 25.17]{Sa04}). 
By Proposition \ref{Ex-EMM}, there exists the EMM $\mb P_1^E$ with respect to $\mb P_1$. Obviously, $\mb P_1^E$ is a martingale measure, and $\mb P_1^E\sim\mb P_1\sim\mb P$, thus $\mb P_1^E\sim\mb P$. This proves the existence of an equivalent martingale measure and hence the claim.\hfill$\Box$
\bg{Cor}
Let $\scr M_{efl}$ be the set of equivalent martingale measures with finite entropy which are preserving the L\'evy property (cf. \eqref{efl-mart} for $(X,\mb F)=(L,\mb F)$). Then, if the market is free of arbitrage, $\scr M_{efl}\not=\emptyset$.
\e{Cor}
\proof The EMM $\mb P_1^E$ with respect to $\mb P_1$ constructed in the proof of Lemma \ref{No-Arb} is equivalent to $\mb P$ and L\'evy preserving. It has also finite entropy $I_T(\mb P_1^E,\mb P)$: 
$$
0\le I_T(\mb P_1^E,\mb P)=I_T(\mb P_1^E,\mb P_1)+\mb E_{\mb P_1^E}\left[\ln\frac{\rmd\mb P_1}{\rmd \mb P}|_{\scr F_T}\right]<+\infty\,,
$$
because of identity \eqref{EM-ent} (for $\mb P_1$ instead of $\mb P$) and the boundedness of $Z_T=\frac{\rmd\mb P_1}{\rmd \mb P}|_{\scr F_T}$ (cf. \eqref{density}).\hfill$\Box$
\bg{Rem} {\rm Some authors (see Esche \& Schweizer \cite{ES05}, Frittelli \cite{Fr00}) use the condition $\scr M_{efl}\not= \emptyset$ as standing hypothesis. We have now seen that this hypothesis is already satisfied if we only assume that the market is free of arbitrage.}
\e{Rem}

\smallskip
\noindent\textit{\textbf{The EMM is always the MEMM}} \ 
As the last preparatory step, now we deal with the easy part of the identification of EMM and MEMM. Let \OFPF\ be a filtered probability space and $(L,\mb F)$ an arbitrary L\'evy process on it. We recall Definition \ref{DEMM} of the EMM $\mb P^E$. 
\bg{Th}\label{EMM=>MEMM}
The EMM $\mb P^E$, if it exists, is the MEMM in the class $\scr M_a$.
\e{Th}
\proof Assume that there exists the EMM $\mathbb{P}^E=\mb P^{\kappa_0}$.
From the above, it is known that $I_T(\mathbb{P}^E, \mathbb{P})=-\log\varphi_T(\kappa_0)<+\infty$.
Therefore, the class $\mathscr{M}_f=\{\mb Q\in\scr M_a: \ I_T(\mb Q,\mb P)<\infty\}$ is not empty.
Let $\mathbb{Q}\in \scr M_a$ be given. Then, it easily follows
\beqas
I_T(\mathbb{Q}, \mathbb{P})&=&\mathbb{E}_{\mb Q}\left[\log\frac{\rmd \mathbb{Q}}{\rmd \mathbb{P}}|_{\scr F_T}\right]\\
&=&\mathbb{E}_{\mb Q}\left[\log\frac{\rmd \mathbb{Q}}{\rmd \mathbb{P}^E}|_{\mathscr{F}_T}+\log\frac{\rmd \mathbb{P}^E}{\rmd \mathbb{P}}|_{\scr F_T}\right]\\
&=&I_T(\mathbb{Q},\mathbb{P}^E)+\mathbb{E}_{\mb Q}[\kappa L_T-\log(\varphi_T (\kappa))]\\ 
&\ge&\kappa\,\mathbb{E}_{\mb Q}[L_T]-\log(\varphi_T (\kappa))=I_T(\mathbb{P}^E, \mathbb{P})\,.
\eeqas
This yields $I_T(\mathbb{P}^E, \mathbb{P})\ge \inf_{\mb Q\in\scr M_a} I_T(\mathbb{Q}, \mathbb{P})\ge I_T(\mathbb{P}^E, \mathbb{P})$ and hence the claim.\hfill$\Box$

\smallskip
The next Corollary is an important observation. 
\bg{Cor}
The EMM $\mb P^E$, if it exists, is the MEMM in the broader class $\tilde{\scr M}_a$.
\e{Cor}
\proof Recalling the definition of $\tilde{\scr M}_a$ in \eqref{abs-mom} and noting that in the proof of Theorem \ref{EMM=>MEMM} it is only used that $\mathbb{E}_{\mb Q}[L_T]=0$ (instead of the martingale property of $(L,\mb F)$ with respect to $\mb Q$), the result follows.\hfill$\Box$

\section{Approximation of the Entropy}\label{sec:App}
Let \OFPF\ be a filtered probabilty space and $(L,\mb F)$ an arbitrary L\'evy process on it admitting $(b,\sig^2,\nu)$ as characteristic triplet.
 
In this section, we prepare the proof of the main result of the paper: If the MEMM $\mb P^*$ exists, then there exists the EMM $\mb P^E$, and both probabilty measures are equal. This is realized by a suitable approximation of the physical probability measure $\mb P$.

We define a sequence $(Y_n)_{n\ge1}$ of Girsanov quantities $Y_n$ by 
\beq\label{GQn}
Y_n(x)=\exp(-\rho_n(x)), \quad x\in{\mathbb{R}}, \ n\ge1,
\eeq
where $\rho_n(x)=x^2/n$ if $|x|>1$ and $0$ otherwise.
We put
\beq\label{density:n}
Z^n_t=\scr E\big((Y_n-1)*(\mu-\mu^p)\big)_t=\exp\lk-\sum_{0<u\le t}\rho_n(\Delta L_u)+t\,\int_\mb R \lk 1-Y_n(x)\rk\,\nu(\rmd x)\rk\,,
\eeq
where $(Y_n-1)*(\mu-\mu^p)$ is the stochastic integral with respect to the compensated Poisson random measure $\mu-\mu^p$ associated with $L$, $\mu$ being the jump measure of $L$ and $\mu^p=\lambda\otimes\nu$ ($\lambda$ Lebesgue measure on $[0,T]$) its predictable compensator, and $\scr E(S)$ denotes the stochastic exponential of a semimartingale $S$. The term on the right hand side is obtained by a straightforward calculation of the stochastic exponential following its definition (cf. \cite[Ch. II.8]{JS00}), observing that $(Y_n-1)*(\mu-\mu^p)=(Y_n-1)*\mu-(Y_n-1)*\mu^p$ because of $Y_n-1\in L^1(\mb R)\cap L^2(\mb R)$ (cf. \cite[Proposition II.1.28]{JS00}). Note that $|\Delta L_u|>1$ for only a finite number of $u\in[0,T]$, hence $\sum_{0<u\le t}\rho_n(\Delta L_u)$ is well defined and finite for all $t\in[0,T]$. In view of $Y_n-1\in L^2(\mb R)$, the stochastic integral $\lk(Y_n-1)*(\mu-\mu^p),\mb F\rk$ is a square integrable martingale and hence its stochastic exponential $\lk Z^n,\mb F\rk$ is a local martingale (cf. \cite[Section II.3]{IkWa} and Proposition \ref{lmm} (i)). From \eqref{density:n}, it follows the estimate
\beq\label{uniformly bounded}
\sup_{0\le t\le T, \ 1\le n}
{Z^n_t}\le\exp{\left(T\int\limits_{\mathbb{R}}{(\exp(-\rho_1(x))-1)\nu(\rmd x)}\right)}< +\infty\,,
\eeq
hence $Z^n$ is uniformly bounded and, in particular, $\lk Z^n,\mb F\rk$ is a martingale with $\mb E_\mb PZ^n_T=1$. Now we define the probabilty measure $\mb P_n$ by $\rmd \mb P_n=Z^n_T\,\rmd \mb P$ which is obviously equivalent to $\mb P$. 
\bg{Prop}\label{P_n:Levy-preserving}
With respect to the probability measure $\mb P_n$, the process $\lk L,\mb F\rk$ is a L\'evy process with characteristic triplet $\lk b,\sig^2,\nu_n\rk$ where the L\'evy measure $\nu_n$ is given by $\nu_n(\rmd x)=\exp(-\rho_n(x))\nu(\rmd x)$.
\e{Prop}
\proof With respect to $\mb P_n$, the process $\lk L,\mb F\rk$ is again a semimartingale, and in view of \cite[Theorem III.3.24]{JS00} the semimartingale characteristics are deterministic and time-independent. Hence $\lk L,\mb F\rk$ is a L\'evy process with respect to $\mb P_n$. The concrete form of the characteristic triplet is given in \cite[III.3.27]{JS00} as claimed above.\hfill$\Box$

\smallskip
By Proposition \ref{P_n:Levy-preserving}, the probability $\mb P_n$ on $\scr F_T$ is L\'evy preserving (see Definition \ref{Lpreserving}), and the L\'evy measure of $(L,\mb F)$ with respect to $\mb P_n$ is $\rmd \nu_n=Y_n(x)\,\rmd\nu(x)$. It follows that $L_T$ has finite exponential moments of arbitrary order with respect to $\mb P_n$ (cf. Sato \cite[Theorem 25.17]{Sa04}). 
Using Proposition \ref{Ex-EMM}, we obtain that there exists the Esscher martingale measure $\mb P_n^E$ with respect to $\mb P_n$. 

For any probability measure $\mb Q$ on $\scr F_T$ such that $I_T(\mb Q,\mb P)<\infty$, we now introduce the finite nonnegative measure $\mb Q_n$ by 
\beq\label{PM-Qn}
\rmd \mathbb{Q}_n=Z^{n}_T\rmd \mathbb{Q}\,.
\eeq 
\begin{Prop}\label{identity1} For any probability measure $\mb Q$ on $\scr F_T$ such that $I_T(\mb Q,\mb P)<\infty$,
$$
I_T(\mathbb{Q},\mathbb{P})=\displaystyle\lim_{n\to\infty} I_T(\mathbb{Q}_n,\mathbb{P}_n)\,.
$$
\end{Prop}
\proof By definition
\beqas
I_T(\mathbb{Q}_n,\mathbb{P}_n)&=&\mb E_{\mb Q_n}\left[\log\frac{\rmd \mathbb{Q}_n}{\rmd \mathbb{P}_n}\right]=\mathbb{E}_{\mb Q_n}\left[\log\frac{\rmd \mathbb{Q}}{\rmd \mathbb{P}}\right]=\mathbb{E}_{\mb Q}\left[\frac{\rmd \mathbb{Q}_n}{\rmd \mathbb{Q}}(\scr F_T)\log\frac{\rmd \mathbb{Q}}{\rmd \mathbb{P}}(\scr F_T)\right]\\
&=&\mathbb{E}_{\mb Q}\left[Z^n_T\log\frac{\rmd \mathbb{Q}}{\rmd \mathbb{P}}(\scr F_T)\right] \ \longrightarrow_{n\to\infty} \ \mathbb{E}_{\mb Q}\left[\log\frac{\rmd \mathbb{Q}}{\rmd \mathbb{P}}(\scr F_T)\right]=I_T(\mathbb{Q}, \mathbb{P})\,,
\eeqas
because of the dominated convergence theorem of Lebesgue: $\log\frac{\rmd \mathbb{Q}}{\rmd \mathbb{P}}(\scr F_T)$ is $\mb Q$-integrable in view of $I_T(\mb Q,\mb P)<\infty$, $\displaystyle\lim_{n\to\infty}{Z^n_T}=1$ and $Z^n_T$ is uniformly bounded by \eqref{density:n} and \eqref{uniformly bounded}.\hfill$\Box$   

\smallskip
In the next proposition, we compare the asymptotic behavior of the entropy of the Esscher measures $\mb P_n^E$ to $\mb P_n$.
\bg{Prop} \label{identity2} It holds
$$
\overline{\lim_{n\to\infty}}I_T(\mathbb{P}_n^E,\mathbb{P})\le\overline{\lim_{n\to\infty}}I_T(\mathbb{P}_n^E,\mathbb{P}_n)\,.
$$
\e{Prop}
\proof An easy calculation shows
\beqas
I_T(\mathbb{P}_n^E,\mathbb{P})&=&\mathbb{E}_{\mb P_n^E}\left[\log\frac{\rmd \mathbb{P}_n^E}{\rmd \mathbb{P}}|_{\scr F_T}\right]\nonumber\\
&=&\mathbb{E}_{P_n^E}\left[\log\frac{\rmd \mathbb{P}_n^E}{\rmd \mathbb{P}_n}|_{\scr F_T}+\log\frac{\rmd \mathbb{P}_n}{\rmd \mathbb{P}}|_{\scr F_T}\right]\\
&=&I_T(\mathbb{P}_n^E,\mathbb{P}_n)+\mathbb{E}_{\mb P_n^E}\left[\log\frac{\rmd \mathbb{P}_n}{\rmd \mathbb{P}}|_{\scr F_T}\right]\\
&=&I_T(\mathbb{P}_n^E,\mathbb{P}_n)+\mathbb{E}_{\mb P_n^E}\left[\log Z^n_T \right].
\eeqas
Here, the linearity of the integral is applied which is possible since $\log\frac{\rmd \mathbb{P}_n^E}{\rmd \mathbb{P}_n}|_{\scr F_T}$ is integrable with respect to $\mb P_n^E$ and the integral of $\log\frac{\rmd \mathbb{P}_n}{\rmd \mathbb{P}}|_{\scr F_T}=\log Z^n_T$ with respect to $\mb P_n^E$ exists and is less than $+\infty$ because $\log Z^n_T$ is bounded from above (see \eqref{density:n}).
Now, by Lebesgue's theorem
$$
\overline{\lim_{n\to\infty}}\mathbb{E}_{P_n^E}\left[\log Z^{n}_T \right]\le\overline{\lim_{n\to\infty}}T\int\limits_{\mathbb{R}}{(1-\exp(-\rho_n(x)))\nu(\rmd x)}= 0,
$$
which yields the claim.\hfill$\Box$

\medskip
Our next goal is to compare the asymptotics of $I_T(\mathbb{Q}_n,\mathbb{P}_n)$ and $I_T(\mathbb{P}_n^E,\mathbb{P}_n)$. Recall that $\mathbb{Q}_n$ is defined by \eqref{PM-Qn}.
\begin{Prop}\label{identity3}
Let $\mathbb{Q}$ be from $\tilde{\mathscr{M}}^{\rm loc}_a$ (cf. \eqref{abs-loc-mom} for $(X,\mb F)=(L,\mb F))$. Then
$$
\overline{\lim_{n\to\infty}}I_T(\mathbb{P}_n^E,\mathbb{P}_n)\le\overline{\lim_{n\to\infty}}I_T(\mathbb{Q}_n,\mathbb{P}_n)\,.
$$
\end{Prop}
\proof By assumption, $\mb Q$ satisfies the local moment condition (see Definition \ref{moment-condition}): There exists an increasing sequence $(\tau_k)$ of $\mb F$-stopping times such that $\{\tau_k=T\}\uparrow_{k\to\infty}\Om$ and $\mb E_\mb Q\left[L_{\tau_k\wedge T}\right]=0$ for all $k\ge1$. The monotonicity of the entropy implies
\beqas
I_T(\mathbb{Q}_n,\mathbb{P}_n)&\ge&I_{\tau_k\wedge T}(\mathbb{Q}_n,\mathbb{P}_n):=\mathbb{E}_{\mb Q_n}\left[\log\frac{\rmd \mathbb{Q}_n}{\rmd \mathbb{P}_n}|_{\scr F_{\tau_k\wedge T}}\right]\\
&=&\mathbb{E}_{\mb Q_n}\left[\log\frac{\rmd \mathbb{Q}_n}{\rmd \mathbb{P}_n^E}|_{\scr F_{\tau_k\wedge T}}+\log\frac{\rmd \mathbb{P}_n^E}{\rmd \mathbb{P}_n}|_{\scr F_{\tau_k\wedge T}}\right]\\
&=&\mathbb{E}_{\mb Q_n}\left[\log\frac{\rmd \mathbb{Q}_n}{\rmd \mathbb{P}_n^E}|_{\scr F_{\tau_k\wedge T}}\right] + \mathbb{E}_{\mb Q_n}\left[\log\frac{\rmd \mathbb{P}_n^E}{\rmd \mathbb{P}_n}|_{\scr F_{\tau_k\wedge T}}\right]\,,
\eeqas
where the linearity of the integral with respect to $\mb Q_n$ is applied which is possible because
\beq\label{entropie:positiv}
\mb E_{\mb Q_n}\left[\log\frac{\rmd \mb Q_n}{\rmd \mb P_n}|_{\scr F_{\tau_k\wedge T}}\right]=I_{\tau_k\wedge T}(\mb Q_n,\mb P^E_n)\ge 0
\eeq
and $\log\frac{\rmd\mb P^E_n}{\rmd\mb P_n}|_{\scr F_{\tau_k\wedge T}}$ is $\mb Q_n$-integrable, as it is easy to verify. Using \eqref{entropy:positiv} and noting that
$$
\log\frac{\rmd \mathbb{P}_n^E}{\rmd \mathbb{P}_n}|_{\scr F_{\tau_k\wedge T}}=-(\tau_k\wedge T)\,\log\mb E_{\mb P_n}\left[\exp\left(\kappa_n\,L_1\right)\right]+\kappa_n\,L_{\tau_k\wedge T}\,,
$$
where $\kappa_n$ is the Esscher parameter of the EMM $\mb P_n^E$ with respect to $\mb P_n$, implies
\beq\label{ineq1}
I_T(\mathbb{Q}_n,\mathbb{P}_n)\ge -\log\mb E_{\mb P_n}\left[\exp\left(\kappa_n\,L_1\right)\right]\;\mb E_{\mb Q_n}\left[\tau_k\wedge T\right]\,+\kappa_n\;\mb E_{\mb Q_n}\left[\,L_{\tau_k\wedge T}\right]\,.
\eeq
Now it follows
$$
\lim_{n\to\infty}\mb E_{\mb Q_n}\left[L_{\tau_k\wedge T}\right]=\lim_{n\to\infty}\mb E_{\mb Q}\left[Z^n_T\,L_{\tau_k\wedge T}\right]=\mb E_{\mb Q}\left[L_{\tau_k\wedge T}\right]=0
$$
in view of Lebesgue's theorem on dominated convergence, the boundedness of $Z^n_T$ and $Z^n_T\longrightarrow_{n\to\infty}1$ (see \eqref{uniformly bounded}) and the choice of $(\tau_k)$. As $(\kappa_n)$ is bounded (see Lemma \ref{entr:conv} below), this yields that the second term of the right-hand side of \eqref{ineq1} converges to zero as $n\to\infty$. Applying \eqref{EM-ent} to $\mb P^E_n$ and $\mb P_n$, the first term on the right hand side of \eqref{ineq1} can be rewritten as $T^{-1}\;I_T(\mathbb{P}_n^E,\mathbb{P}_n)\;\mb E_{\mb Q_n}\left[\tau_k\wedge T\right]$.
Hence
$$
\overline{\lim_{n\to\infty}}I_T(\mathbb{Q}_n,\mathbb{P}_n)\ge
T^{-1}\;\overline{\lim_{n\to\infty}}I_T(\mathbb{P}_n^E,\mathbb{P}_n)\;\mb E_{\mb Q}\left[\tau_k\wedge T\right]\,.
$$
Letting $k$ converge to infinity, yields the claim. \hfill$\Box$

\medskip
Before we come to the next lemma, we introduce the moment generating functions $\ph_{T,n}$ with respect to the probability measures $\mb P_n$ (see Appendix): 
\beq\label{moment generating:n}
\ph_{T,n}(\kappa)=\mb E_{P_n}\left[\exp\left(\kappa\,L_T\right)\right], \quad\kappa\in\mb R,\quad n\ge1\,.
\eeq
\bg{Le}\label{entr:conv} Let $\kappa_n$ be the Esscher parameter of the EMM $\mb P_n^E$ with respect to $\mb P_n$. Then it follows:

{\rm (i)} \ The set $\{\kappa_n: \ n\ge 1\}$ is bounded.

{\rm (ii)} \ $\kappa_0:=\lim_{n\to\infty} \kappa_n$ exists and $\lim_{n\to\infty} \ph_{T,n}(\kappa_n)=\ph_T(\kappa_0)$. Moreover, $\ph_T$ reaches its minimum at $\kappa_0$. 
\e{Le}
\proof By definition of the EMM $\mathbb{P}_n^E$, we have:
$$
\mathbb{E}_{\mb P_n}\left[L_T\exp\left(\kappa_n L_T\right)\right]=0\,.
$$
The function $\varphi_{T,n}$ is twice continuously differentiable and 
$$
\varphi_{T,n}^{\prime}(\kappa)=\mathbb{E}_{\mb P_n}\left[L_T\exp\left(\kappa L_T\right)\right],\quad \varphi_{T,n}^{\prime\prime}(\kappa)=\mathbb{E}_{\mb P_n}\left[L^2_T\exp\left(\kappa L_T\right)\right],\quad \kappa\in\mathbb{R}
$$
(see Propositions \ref{prop:phi} and \ref{prop:psi}). From the definition of $\kappa_n$ now follows that $\varphi_{T,n}^{\prime}(\kappa_n)=0$ and, because of $\varphi_{T,n}^{\prime\prime}(\kappa_n)>0$, $\varphi_{T,n}$ reaches its minimum at the point $\kappa_n$. Define sets $K_n$ by
$$
K_n:=\left\{\kappa\in\mathbb{R}: \, \varphi_{T,n}(\kappa)\exp\left(-T\int_{\mathbb{R}}{(1-\exp(-\rho_n(x)))\nu(\rmd x)}\right)\le1\right\}
$$
and notice that 
\begin{eqnarray*}
\lefteqn{\varphi_{T,n}(\kappa)\exp\left(-T\int_{\mathbb{R}}{(1-\exp(-\rho_n(x)))\nu(\rmd x)}\right)}\\
&=&\mb E_{\mb P_n}\left[\exp\left(\kappa L_T\right)\right]\, \exp\left(-T\int_{\mathbb{R}}{(1-\exp(-\rho_n(x)))\nu(\rmd x)}\right)\\
&=&\mb E_{\mb P}[Z_T^{n}\exp(\kappa L_T)] \, \exp\left(-T\int_{\mathbb{R}}{(1-\exp(-\rho_n(x)))\nu(\rmd x)}\right)\\
&=&\mathbb{E}_{\mb P}\left[\exp\left(-\sum_{0<u\le T}\rho_n(\Delta L_u)\right)\exp(\kappa L_T)\right],
\end{eqnarray*}
where the representation \eqref{density:n} for the density process $Z^n$ is used. The last term being monotonically increasing in $n$, shows that the sets $K_n$ are monotonically decreasing and, in particular, that $K_n\subseteq K_1$ for all $n\ge 1$.
The function $\varphi_{T,n}$ is strictly positive, $\varphi_{T,n}(0)=1$, hence its minimum is not larger than 1 and therefore $\kappa_n\in K_n\subseteq K_1$, for all $n\ge1$. It can easily be verified that $K_1$ is compact: $K_1$ is bounded because $\lim_{|\kappa|\to\infty}\ph_{T,n}(\kappa)=+\infty$ (see Proposition \ref{prop:phi}); $K_1$ is closed because $\varphi_{T,n}$ is continuous (see Proposition \ref{prop:phi}). This proves (i). For proving (ii), let $(\kappa_{n'})$ be a converging subsequence of $(\kappa_{n})$ with limit $\kappa_0$. Using Fatou's lemma and then Lebesgue's dominated convergence theorem, yields
\beqas
\ph_T(\kappa_0)&=&\mb E_\mb P\left[\kappa_0\,L_T\right]\le\liminf_{n'\to\infty} \mb E_\mb P\left[Z^{n'}_T\,\exp\left(\kappa_{n'}\,L_T\right)\right]\\
&=&\liminf_{n'\to\infty}\ph_{T,n'}(\kappa_{n'})\le\limsup_{n'\to\infty}\ph_{T,n'}(\kappa_{n'})\\
&\le& \limsup_{n'\to\infty}\ph_{T,n'}(\kappa)=\ph_T(\kappa)
\eeqas
for any $\kappa\in\mb R$ such that $\ph_T(\kappa)<+\infty$. This implies
\beq\label{chain}
\inf_{\kappa\in\mb R} \ph_T(\kappa)\le\ph_T(\kappa_0)\le\lim_{n'\to\infty}\ph_{T,n'}(\kappa_{n'})\le\inf_{\kappa\in\mb R} \ph_T(\kappa)\,.
\eeq
Hence $\ph_T$ reaches its minimum in $\kappa_0$. Similarly, for another subsequence $(\kappa_{n''})$ of $(\kappa_{n})$ converging to $\kappa_1$, $\ph_T$ reaches its minimum in $\kappa_1$ and, hence, $\ph_T(\kappa_0)=\ph_T(\kappa_1)$. Since $\ph_T$ is strictly convex on the interval on which $\ph_T$ takes finite values, it follows that $\kappa_0=\kappa_1$. This yields that $(\kappa_n)$ converges to the unique minimal point $\kappa_0$ of $\ph_T$, and the above chain of inequalities \eqref{chain}, now for the sequence $(\kappa_n)$, shows that $\ph_T(\kappa_0)=\lim_{n\to\infty}\ph_{T,n}(\kappa_n)$. This proves (ii).\hfill$\Box$
\bg{Rem} {\rm For the approximation procedure, we have choosen the sequence of functions $\rho_n$ for the Girsanov quantities $Y_n$ (see \eqref{GQn}) in a particularly simple way. However, for the approximation of the minimal entropy in special situations important for applications, it could be of interest to make another choice. Here are the conditions on $(\rho_n)$ which are needed that the approximation works well:

\medskip
(1) \ $0\leq \rho_{n+1}$, $\rho_{n+1}\le \rho_n$, $\displaystyle\lim_{n\to\infty}\rho_n(x)=0$;

\smallskip
(2) \ $\displaystyle\frac{ |x|}{\rho_n(x)}\to 0$ if $|x|\to \infty$, i.e., $|x|={\rm o}(\rho_n(x))$ if $|x|\to \infty$, for all $n\ge 1$;

\smallskip
(3) \ $1-\exp(-\rho_n)\in L^1(\nu)$.

\medskip
\noindent The proofs are similar and left to the reader.
}
\e{Rem}

\section{Conclusions}\label{sec:con}
By $\kappa_0$ we denote the unique point in which $\ph_T$ attains its minimum (see Proposition \ref{prop:phi}).
\begin{Th}\label{identity4}
The following identity holds:
$$
\lim_{n\to\infty}I_T(\mathbb{P}^E_n, \mathbb{P})=\lim_{n\to\infty}I_T(\mathbb{P}^E_n, \mathbb{P}_n)=\inf_{\mb Q\in\tilde{\mathscr{M}}^{\rm loc}_a} I_T(\mathbb{Q}, \mathbb{P})=-\log \ph_T(\kappa_0)\,.
$$
\end{Th}
\proof For any $\mathbb{Q}\in\tilde{\mathscr{M}}_a^{\rm loc}$, from Propositions~\ref{identity3} and ~\ref{identity1}, it follows
$$
\underline{\lim}_{n\to\infty}I_T(\mathbb{P}^E_n, \mathbb{P})\le\overline{\lim_{n\to\infty}}I_T(\mathbb{P}^E_n, \mathbb{P})\le
\overline{\lim_{n\to\infty}}I_T(\mathbb{Q}_n, \mathbb{P}_n)=I_T(\mathbb{Q}, \mathbb{P})\,.
$$
Hence
$$
\inf_{\mb Q\in\tilde{\mathscr{M}}_a^{\rm loc}}I_T (\mathbb{Q}, \mathbb{P})\le\underline{\lim}_{n\to\infty}I_T(\mathbb{P}^E_n, \mathbb{P})\le
\overline{\lim_{n\to\infty}}I_T(\mathbb{P}^E_n, \mathbb{P})\le\inf_{\mb Q\in\tilde{\mathscr{M}}_a^{\rm loc}}I_T (\mathbb{Q}, \mathbb{P})\,,
$$
proving
$$
\lim_{n\to\infty}I_T(\mathbb{P}^E_n, \mathbb{P})=\inf_{\mb Q\in\tilde{\mathscr{M}}_a^{\rm loc}}I_T(\mathbb{Q}, \mathbb{P})\,.
$$
In a similar way, it can be shown that
$$
\overline{\lim_{n\to\infty}}I_T(\mathbb{P}^E_n, \mathbb{P}_n)=\inf_{\mb Q\in\tilde{\mathscr{M}}_a^{\rm loc}}I_T(\mathbb{Q}, \mathbb{P})\,.
$$
Now Lemma \ref{entr:conv} (ii) implies that the limit of $I_T(\mathbb{P}^E_n, \mathbb{P}_n)=-\log\ph_{T,n}(\kappa_n)$ exists and is equal to $-\log \ph_T(\kappa_0)$. This  proves the theorem. \hfill$\Box$
\bg{Cor}\label{identities}
$$
\inf_{\mb Q\in\scr M_{efl}} I_T(\mb Q,\mb P)
=\inf_{\mb Q\in\mathscr{M}_a}I_T(\mathbb{Q}, \mathbb{P})
=\inf_{\mb Q\in\mathscr{M}^{\rm loc}_a}I_T(\mathbb{Q}, \mathbb{P})
=\inf_{\mb Q\in\tilde{\mathscr{M}}^{\rm loc}_a}I_T(\mathbb{Q}, \mathbb{P})\,.
$$
\e{Cor}
\proof Noting that the Esscher martingale measures $\mb P^E_n$ belong to $\scr M_{efl}$, this follows from Theorem \ref{identity4}.\hfill$\Box$
\bg{Cor}[Preservation of L\'evy Property]
The class $\mathscr{M}_{efl}$ of equivalent martingale measures $\mb Q$ with finite entropy which preserve the L\'evy property is a sufficient subclass for determining the minimal entropy:
$$
\inf_{\mb Q\in\scr M_{efl}} I_T(\mb Q,\mb P)=\inf_{\mb Q\in\tilde{\mathscr{M}}^{\rm loc}_a}I_T(\mathbb{Q}, \mathbb{P})\,.
$$ 
\e{Cor}

Essentially, this was proven by Esche \& Schweizer \cite{ES05} by averaging arbitrary Girsanov parameters but it was not stated explicitly. In fact, from this they derived that the MEMM, if it exixts, preserves the L\'evy property (see \cite[Theorem A]{ES05}).
\bg{Cor}[Moment Problem]\label{Moment Problem}
Minimizing the entropy over the class $\scr M_a$ (resp., $\scr M_a^{\rm loc})$)  is equivalent to minimizing the entropy over the broader class $\tilde{\scr M}_a$ (resp., $\tilde{\scr M}^{\rm loc}_a$):
$$
\inf_{\mb Q\in\mathscr{M}_a}I_T(\mathbb{Q}, \mathbb{P})=\inf_{\mb Q\in\tilde{\mathscr{M}}_a}I_T(\mathbb{Q}, \mathbb{P})\,.
$$
(resp., 
$$
\inf_{\mb Q\in\mathscr{M}^{\rm loc}_a}I_T(\mathbb{Q}, \mathbb{P})=\inf_{\mb Q\in\tilde{\mathscr{M}}^{\rm loc}_a}I_T(\mathbb{Q}, \mathbb{P}))\,.
$$
\e{Cor}
\bg{Rem}{\rm
Searching for the minimal entropy $\inf_{\mb Q\in\tilde{\mathscr{M}}_a}I_T(\mathbb{Q}, \mathbb{P})$ for the L\'evy market $(L,\mb F)$ is actually a problem in the one-step financial market $(X,\tilde{\mb F})$ defined by $X_0=0$, $\tilde{\scr F}_0=\{\emptyset,\Om\}$; $X_1=L_T$, $\tilde{\scr F}_1=\scr F_T$ (moment problem). Indeed, $\tilde{\mathscr{M}}_a$ is just the class of absolutely continuous probability measures $\mb Q$ rendering $(X,\tilde{\mb F})$ a $\mb Q$-martingale. Thus, the meaning of Corollary \ref{Moment Problem} is that the problem of finding the minimal entropy in the class $\mathscr{M}_a$ for $(L,\mb F)$ is equivalent to the problem of finding the minimal entropy in the class $\tilde{\mathscr{M}}_a$ for the one-step market $(X,\tilde{\mb F})$. As we shall see in Corollary \ref{id-MEMM}, if the MEMM for one problem exists, then the MEMM also exists for the other and both coincide, and coincide with the EMM, too. Note that the EMM for the L\'evy market $(L,\mb F)$ and the EMM for the one-step market $(X,\tilde{\mb F})$ are the same. 
}  
\e{Rem}

Now we come to the main theorem asserting that, for an arbitrary linear L\'evy market $(L,\mb F)$ on a filtered probability space \OFPF\ the minimal entropy martingale measure MEMM coincides with the Esscher martingale measure EMM. For the definition of the classes $\hat{\scr M}^{\rm loc}_a$ and $\tilde{\scr M}^{\rm loc,0}_a$, see \eqref{abs-loc-mom-int} and \eqref{abs-loc-hat} for $(X,\mb F)=(L,\mb F)$.
\bg{Th}\label{mainth}
The MEMM $\mb P^*$ in the class $\hat{\scr M}^{\rm loc}_a$ coincides with the EMM $\mb P^E$: $\mb P^*$ exists if and only if $\mb P^E$ exists, and in this case $\mb P^*=\mb P^E$.  
\e{Th}
\proof If The EMM $\mb P^E$ exists, then it is the MEMM $\mb P^*$ in ${\scr M}_a$ by Theorem \ref{EMM=>MEMM}. Corollary \ref{identities} implies that $\mb P^*$ is the MEMM in $\tilde{\scr M}^{\rm loc}_a$ and hence also in $\hat{\scr M}^{\rm loc}_a$. Conversely, let $\mb P^*$ be the MEMM in $\hat{\scr M}^{\rm loc}_a$. Let $\kappa_0$ be the unique point in which $\ph_T$ reaches its minimum and denote by $\mb P^{\kappa_0}$ the Esscher measure with Esscher parameter $\kappa_0$ which is well defined because $\mb E_{\mb P}\left[\exp\left(\kappa_0\,L_T\right)\right]<+\infty$. The first step now consists in verifying that $I_T(\mb P^*,\mb P^{\kappa_0})=0$: \ Since $\mb P^*$ belongs to $\tilde{\mathscr{M}}^{\rm loc}_a$, there exists an increasing sequence $(\tau_k)$ of $\mb F$-stopping times such that $\{\tau_k=T\}\uparrow_{k\to\infty}\Om$ and $\mb E_{\mb P^*}\left[L_{\tau_k\wedge T}\right]=0$, and it follows 
\beqas
I_T(\mathbb{P}^*, \mathbb{P})&\ge&I_{\tau_k\wedge T}(\mathbb{P}^*, \mathbb{P})=\mathbb{E}_{\mb P^*}\left[\log\frac{\rmd \mathbb{P}^*}{\rmd \mathbb{P}}|_{\mathscr{F}_{\tau_k\wedge T}}\right]\\
&=&\mathbb{E}_{\mb P^*}\left[\log\frac{\rmd \mathbb{P}^*}{\rmd \mathbb{P}^{\kappa_0}}|_{\mathscr{F}_{\tau_k\wedge T}}\right] + \mathbb{E}_{\mb P^*}\left[\log\frac{\rmd \mathbb{P}^{\kappa_0}}{\rmd \mathbb{P}}|_{\mathscr{F}_{\tau_k\wedge T}}\right]\\
&=&\mathbb{E}_{\mb P^{\kappa_0}}\left[\frac{\rmd \mathbb{P}^*}{\rmd \mathbb{P}^{\kappa_0}}|_{\mathscr{F}_{\tau_k\wedge T}}\,\log\frac{\rmd \mathbb{P}^*}{\rmd \mathbb{P}^{\kappa_0}}|_{\mathscr{F}_{\tau_k\wedge T}}\right]+\mathbb{E}_{\mb P^*}\left[\kappa_0 L_{\tau_k\wedge T}-(\tau_k\wedge T)\,\log\varphi_1(\kappa_0)\right]\\
&=&\mathbb{E}_{\mb P^{\kappa_0}}\left[\frac{\rmd \mathbb{P}^*}{\rmd \mathbb{P}^{\kappa_0}}|_{\mathscr{F}_{\tau_k\wedge T}}\,\log\frac{\rmd \mathbb{P}^*}{\rmd \mathbb{P}^{\kappa_0}}|_{\mathscr{F}_{\tau_k\wedge T}}\right]-\mathbb{E}_{\mb P^*}[\tau_k\wedge T]\,\log\varphi_1(\kappa_0)\,.
\eeqas
Note that the function $x\mapsto x\log x$ is bounded from below, and passing to the limit for $k\to\infty$ and using Fatou's lemma, yields
\beqas
I_T(\mathbb{P}^*, \mathbb{P})&\ge& \underline{\lim}_{k\to\infty}\mathbb{E}_{\mb P^{\kappa_0}}\left[\frac{\rmd \mathbb{P}^*}{\rmd \mathbb{P}^{\kappa_0}}|_{\mathscr{F}_{\tau_k\wedge T}}\,\log\frac{\rmd \mathbb{P}^*}{\rmd \mathbb{P}^{\kappa_0}}|_{\mathscr{F}_{\tau_k\wedge T}}\right]-T\,\log\varphi_1(\kappa_0)\\
&\ge& \mathbb{E}_{\mb P^{\kappa_0}}\left[\frac{\rmd \mathbb{P}^*}{\rmd \mathbb{P}^{\kappa_0}}|_{\mathscr{F}_{\tau_k\wedge T}}\,\log\frac{\rmd \mathbb{P}^*}{\rmd \mathbb{P}^{\kappa_0}}|_{\mathscr{F}_{\tau_k\wedge T}}\right]-T\,\log\varphi_1(\kappa_0)\\
&=&I_T(\mathbb{P}^*, \mathbb{P}^{\kappa_0})-\log\varphi_T(\kappa_0)\,.\\
\eeqas
By assumption, $\mb P^*$ is the MEMM in $\hat{\scr M}_a^{\rm loc}$. Since ${\scr M}_a^{\rm loc}\subseteq\hat{\scr M}_a^{\rm loc}\subseteq\tilde{\scr M}_a^{\rm loc}$, Theorem \ref{identity4} and Corollary \ref{identities} imply $I_T(\mathbb{P}^*, \mathbb{P})=-\log\varphi_T(\kappa_0)$, consequently,
$$
I_T(\mathbb{P}^*, \mathbb{P})\ge I_T(\mathbb{P}^*, \mathbb{P}^{\kappa_0})+I_T(\mathbb{P}^*, \mathbb{P})\,,
$$
and hence $I_T(\mathbb{P}^*, \mathbb{P}^{\kappa_0})=0$, finishing the proof of the first step. The entropy is zero just in case when the measures coincide on the given $\sigma$-algebra. Consequently, $\mb P^*=\mb P^{\kappa_0}$. This means that the Esscher measure $\mb P^{\kappa_0}$, defined by its parameter $\kappa_0$, coincides with the MEMM $\mathbb{P}^*$. To finish the proof, it is sufficient to verify that $\mb P^{\kappa_0}$ is in fact a martingale measure, namely, the Esscher martingale measure with parameter $\kappa_0$.

Indeed, by assumption, the probability measure $\mb P^*=\mb P^{\kappa_0}$ belongs to $\hat{\scr M}_a^{\rm loc}:=\tilde{\scr M}_a^{\rm loc,0}\cup\scr M_a^{\rm loc}$. Now, if $\mb P^*=\mb P^{\kappa_0}$ belongs to $\tilde{\scr M}_a^{\rm loc,0}$, the Esscher measure $\mb P^{\kappa_0}$ satisfies the assumptions of Lemma \ref{LocMoment} (ii). Hence $(L,\mb F)$ is a martingale with respect to $\mb P^{\kappa_0}$.  If, however, $\mb P^*=\mb P^{\kappa_0}$ belongs to $\scr M_a^{\rm loc}$, $(L,\mb F)$ is a local martingale with respect to $\mb P^{\kappa_0}$ and, since $\mb P^{\kappa_0}$ is L\'evy preserving, it is actually a martingale by Lemma \ref{LocMoment} (i). This concludes the proof of the theorem.\hfill$\Box$
\bg{Cor}\label{id-MEMM}
If the MEMM $\mb P^*$ in $\hat{\scr M}^{\rm loc}_a$ exists, then it exists in any of the classes $\tilde{\scr M}_a$, ${\scr M}^{\rm loc}_a$, ${\scr M}_a$ and ${\scr M}_{efl}$, and is the same.
\e{Cor}
This is obviously satisfied, since $\mb P^*=\mb P^E\in{\scr M}_{efl}$ and the latter is a subclass of all the classes of probability measures appearing in the formulation of the corollary.

\bg{Cor}\label{mainthCor2}
Suppose that $\ph_T$ is not identically equal to $+\infty$ on $\mb R\setminus\{0\}$. Then, the MEMM $\mb P^*$ in the (broader) class $\tilde{\scr M}^{\rm loc}_a$ coincides with the EMM $\mb P^E$.
\e{Cor}
\proof The proof can be given along similar lines as for Theorem \ref{mainth}: \ If $\mb P^*$ is the MEMM in $\tilde{\scr M}^{\rm loc}_a$, then it can be shown as above that $\mb P^*=\mb P^{\kappa_0}$. Now it suffices to verify that $\mb P^{\kappa_0}$ is the EMM. For this, some properties of the functions $\ph_T$ and $\psi_T$ are needed (cf. Appendix). Recall that $I:=\{\kappa\in\mb R: \ \ph_T(\kappa)<+\infty\}$ is an interval with endpoints denoted by $a, b$ with $\infty\le a\le b\le+\infty$. In view of the assumption that $\ph_T$ is not identically equal to $+\infty$ on $\mb R\setminus\{0\}$, it follows $a<b$. If $\kappa_0$ is an inner point of $I$, since $\ph_T$ attains its minimum at $\kappa_0$, then $\psi_T(\kappa_0)=\ph_T^\prime(\kappa_0)=0$, meaning that 
$$
\mb E_{\mb P^{\kappa_0}}\left[L_T\right]=\frac{\mb E_\mb P\left[L_T\,\exp\left(\kappa_0\, L_T\right)\right]}{\mb E_\mb P\left[\exp\left(\kappa_0\, L_T\right)\right]}=\frac{\psi_T(\kappa_0)}{\ph_T(\kappa_0)}=0\,.
$$
Hence $\mb P^{\kappa_0}$ is the EMM. If $a<\kappa_0=b<+\infty$, then $-\infty<\psi_T(\kappa)=\ph_T^\prime(\kappa)<0$ for every $\kappa$ such that $a<\kappa<\kappa_0$ which yields $-\infty<\psi_T(\kappa_0)\le 0$. Because of the definition of $\psi_T$, this means that $L_T$ is integrable with respect to $\mb P^{\kappa_0}$. Moreover, in view of $\mb P^*=\mb P^{\kappa_0}$ and the assumption $\mb P^*\in\tilde{\scr M}^{\rm loc}_a$, it follows that $\mb P^{\kappa_0}\in\tilde{\scr M}_a^{\rm loc,0}$. Now, since $(L,\mb F)$ is a L\'evy process with respect to $\mb P^{\kappa_0}$, Lemma \ref{LocMoment} (ii) implies that $(L,\mb F)$ is a martingale with respect to $\mb P^{\kappa_0}$. Thus, $\mb P^{\kappa_0}$ is the EMM. Similarly, the case $-\infty<a=\kappa_0<b$ can be treated. This shows that $\mb P^*=\mb P^{\kappa_0}$ is the EMM.\hfill$\Box$.
\bg{Cor}\label{mainthCor3}
Suppose that $I:=\{\kappa\in\mb R: \ \ph_T(\kappa)<+\infty\}=\{0\}$. Then the following conditions are equivalent:

\smallskip
{\rm (i)} The MEMM $\mb P^*$ in $\hat{\scr M}^{\rm loc}_a$ exists and is equal to the EMM $\mb P^E$.

\smallskip
{\rm (ii)} $\mb P$ is a martingale measure for $(L,\mb F)$.

\smallskip
{\rm (iii)} $\mb E_\mb P\left[L_t\right]=0$ for $t=T$, and hence for all $t\in[0,T]$.

\smallskip
{\rm (iv)} The EMM $\mb P^E$ exists and is equal to $\mb P$. 

\smallskip
\noindent If one, and therefore all, of these conditions {\rm (i)--(iv)} is satisfied, then $\mb P$ is the MEMM in $\hat{\scr M}^{\rm loc}_a$. Otherwise, the MEMM in the class $\hat{\scr M}^{\rm loc}_a$ and the EMM do not exist.
\e{Cor}  
\bg{Ex}\label{Ex1}{\rm
On a probability space \OFP, let $(L,\mb F)$ be a symmetric $\alpha$-stable process with parameter $\alpha$: \ $0<\alpha<2$. With respect to the standard truncation function, the characteristic triplet is $(0,0,\nu_\alpha)$ where $\nu_\alpha(\rmd x)=|x|^{-\alpha-1}\,\rmd x$. Clearly, $\mb E_\mb P\left[\exp\left(\kappa\,L_T\right)\right]=+\infty$ for all $\kappa\not=0$, so the assumption of Corollary \ref{mainthCor3} is satisfied.

(1) Consider the case $0<\alpha\le 1$. Note that $(L,\mb F)$ is a Cauchy process if $\alpha=1$. Then $\mb E_\mb P\left[L^+_T\right]=\mb E_\mb P\left[L^-_T\right]=+\infty$, hence $\mb E_\mb P\left[L_T\right]$ does not exist. In view of Corollary \ref{mainthCor3}, the MEMM $\mb P^*$ and the EMM $\mb P^E$ do not exist. The function $\ph_T$ is reaching its minimum at $\kappa_0=0$, and $-\log\ph_T(\kappa_0)=0$. The Esscher measure $\mb P^{\kappa_0}=\mb P$ is not a martingale measure, because of the absence of the expectation.  However, by Theorem \ref{identities}, there exists a sequence $(\mb P^E_n)$ of L\'evy preserving equivalent martingale measures such that $\lim_{n\to\infty}I_T(P^E_n,\mb P)=0$ and one is tempted to say that ``$(\mb P^E_n)$ converges to $\mb P$ in entropy". The physical probability measure $\mb P$ is not the MEMM (and not the EMM) only because it misses the existence of the expectation of $L_T$. 

(2) Now consider the case $1<\alpha<2$. Then $\mb E_\mb P\left[|L_T|\right]<+\infty$ and $\mb E_\mb P\left[L_T\right]=0$. Hence $(L,\mb F)$ is a martingale with respect to $\mb P$. Setting $\kappa_0=0$, yields that $\mb P=\mb P^{\kappa_0}=\mb P^E$ is the EMM and hence the MEMM in $\tilde{\mathscr{M}}^{\rm loc}_a$.

}

\e{Ex}

\begin{appendix}
\appendix
\renewcommand{\thesection}{\Alph{section}}
\setcounter{equation}{0}
\stepcounter{section}
\renewcommand{\theequation}{\thesection.\arabic{equation}}
\section*{Appendix A: \ Moment Generating Functions}
\label{sec:app}
\small{
Let $(\Omega, \scr F, \mathbb{P})$ be a probability space and  $\xi$ a random variable defined on it. In the following, it will always be assumed that $\mb P(\{\xi>0\})>0$ and $\mb P(\{\xi<0\})>0$. Note that, in case of a L\'evy process $(L,\mb F)$ on $[0,T]$ and $\xi=L_T$, this corresponds to the no-arbitrage condition of the L\'evy market.   

Define the moment generating function $\ph$ by 
\beq\label{def.phi.app}
\varphi(\kappa)=\mathbb{E}\left[\exp\left(\kappa\,\xi\right)\right], \quad \kappa\in\mathbb{R},
\eeq
and let $I=\left\{\kappa\in\mathbb{R}: \, \varphi(\kappa)<+\infty\right\}$. Then $\ph$ is strictly convex on $I$ and $0\in I$. Therefore, $I$ is an interval with endpoints $a$ and $b$, $-\infty\le a\le 0\le b\le+\infty$. The set of interior points of $I$ is denoted by $I^0$. Clearly, $I^0=(a,b)$. Note that $I^0$ can be the empty set which happens in the case if $a=b=0$. 
\begin{Prop}\label{prop:phi} {\rm (i)} \ $\lim_{|k|\to\infty}\ph(\kappa)=+\infty$;

\smallskip
{\rm (ii)} \ $\varphi$ is continuous on $I$;

\smallskip
{\rm (iii)} \ $\varphi$ is infinitely often differentiable on $(a,b)$; 

\smallskip
{\rm (iv)} \ $\varphi$ has a unique minimum point $\kappa_0$.
\end{Prop}
\proof (i) easily follows from the no-arbitrage condition on $\xi$. (ii) and (iii) are applications of Lebegue's dominated convergence theorem. For showing (iv), let $K=\{\kappa\in\mb R: \ \ph(\kappa)\le1\}$. By (i) and (ii), $K$ is bounded and closed. Hence, $K$ is compact. The function $\ph$ being continuous on $K$, yields that there exists a minimum point $\kappa_0$ on $K$ which is also a minimum point on $\mb R$. Uniqueness follows from the strict convexity of $\ph$ on $I$.\hfill$\Box$  

\smallskip
For given $\kappa\in\mb R$, introduce the notation 
\beq\label{function:psi}
\psi(\kappa):=\mathbb{E}\left[\xi\exp\left(\kappa\xi\right)\right]\,,
\eeq
provided that the integral on the right-hand side exists (it may be equal to $+\infty$ or $-\infty$). Define $E:=\{\kappa\in\mb R: \ -\infty<\psi(\kappa)<+\infty\}$. Obviously, $I^0\subseteq E\subseteq I$. Note that the function $\kappa\mapsto\xi\,\exp\left(\kappa\xi\right)$ is increasing. As a consequence, if $E\not=\emptyset$, it is easy to verify that the integrals \eqref{function:psi} exist for all $\kappa\in\mb R$ and are monotonically increasing in $\kappa$.
\bg{Prop}\label{prop:psi} Suppose that $E\not=\emptyset$.

\smallskip
{\rm (i)} $\psi$ is monotonically increasing, $\lim_{\kappa\to -\infty}\psi(\kappa)=-\infty$, and $\lim_{\kappa\to +\infty}\psi(\kappa)=+\infty$.

\smallskip
{\rm (ii)} \ $\psi$ is continuous on $[a,b]$ (in the topology of the extended real line).

\smallskip
{\rm (iii)} \ $\ph^\prime(\kappa)=\psi(\kappa), \ \kappa\in (a,b)$.

\smallskip
{\rm (iv)} \ If $\psi(b)<+\infty$ (resp., $-\infty<\psi(a)$), then $\ph$ is differentiable at $b$ from the left (resp., at $a$ from the right) and $\ph^\prime(b)=\psi(b)$ (resp., $\ph^\prime(a)=\psi(a))$.  

\smallskip
{\rm (v)} \ $\psi$ is infinitely often differentiable on $(a,b)$. In particular, 
\beq\label{second derivative}
\ph^{\prime\prime}(\kappa)=\psi^\prime(\kappa)=\mb E\left[\xi^2\exp\left(\kappa\xi\right)\right]>0,\quad\mbox{for all } \ \kappa\in (a,b)\,.
\eeq
\e{Prop}
\proof The first part of (i) follows from the existence of the integrals \eqref{function:psi} for $\kappa\in E\not=\emptyset$ and the monotonicity of the function $\kappa\mapsto\xi\exp\left(\kappa\,\xi\right)$. For the second part of (i), note that the integral of the positive part of $\xi\,\exp\left(\kappa\,\xi\right)$,
$$
E\left[\xi^+\exp\left(\kappa\xi\right)\right] =E\left[\xi^+\I_{\{\xi>0\}}\exp\left(\kappa\xi^+\right)\right]\,,
$$
converges to $+\infty$ as $\kappa\to+\infty$ in view of monotone convergence as well as the integral of the negative part,
$$
E\left[\xi^-\exp\left(\kappa\,\xi\right)\right] =E\left[\xi^-\I_{\{\xi<0\}}\exp\left(\kappa\left(-\xi^-\right)\right)\right]\,,
$$
converges to zero as $\kappa\to+\infty$ by Lebesgue's dominated  convergence theorem in view of the assumption that $E\not=\emptyset$. For  $\kappa\to-\infty$, the statement is verified analogously. The proofs of (ii) -- (v) are in a similar fashion, repeatedly using Lebegue's dominated convergence theorem. Details are omitted.\hfill$\Box$

\smallskip
The next proposition characterizes the unique minimal point $\kappa_0\in\mb R$ of $\ph$ in terms of the function $\psi$.
\bg{Prop}\label{char:min}
$\kappa_0\in\mb R$ is the unique minimal point of $\ph$ if and only if one of the following mutually exclusive cases occurs:

\smallskip
{\rm (i)} \ $a=b=0$ and $\kappa_0=0$.

\smallskip
{\rm (ii)} \ $a<b$, $\kappa_0\in[a,b]\cap \mb R$, and $\psi(\kappa_0)=0$.

\smallskip
{\rm (iii)} \ $a<b<+\infty$, $\kappa_0=b$, and $\psi(b)<0$.

\smallskip
{\rm (iv)} \ $-\infty<a<b$, $\kappa_0=a$, and $0<\psi(a)$.
\e{Prop}
\proof Using Propositions \ref{prop:phi} and \ref{function:psi} and, in particular, \eqref{second derivative}, the result follows by elementary calculus.\hfill$\Box$

\smallskip
Points $\kappa_0\in[a,b]\cap\mb R$ such that $\kappa_0\in E\subseteq I$ and $\psi(\kappa_0)=0$ are of special interest. We call them Esscher parameters. This comes from the fact that they determine the Esscher martingale transform $\rmd\mb P^{\kappa_0}=Z_1\rmd\mb P$ with $Z_1:=\exp(\kappa_0\xi)/\mb E\left[\exp(\kappa_0\xi)\right]$ for the one-step model $(X,\mb F)$, $X_0=0$, $\scr F_0=\{\emptyset,\Om\}$; $X_1=\xi$, $\scr F_1=\scr F$. The following proposition gives necessary and sufficient conditions for the existence of an Esscher parameter in terms of the set $E$ and the function $\psi$. Note that an Esscher parameter $\kappa_0$ is always the minimum point of $\ph$  and is therefore unique (see Proposition \ref{prop:phi}). However, the converse is not true. This follows from the characterization of the minimal point of $\ph$ given in Proposition \ref{char:min}.
\bg{Prop}\label{Esscher Parameter}
The Esscher parameter $\kappa_0$ exists in the following cases:

\smallskip
{\rm (i)} $-\infty\le a<b\le+\infty$ and $E=(a,b)$.

\smallskip
{\rm (ii)} $-\infty\le a<b<+\infty$ and $E=(a,b]$ with $\psi(b)\ge0$.

\smallskip
{\rm (iii)} $-\infty<a<b\le+\infty$ and $E=[a,b)$ with $\psi(a)\le0$.

\smallskip
{\rm (iv)} $-\infty<a\le b<+\infty$ and $E=[a,b]$ with $\psi(a)\le0\le\psi(b)$.

\smallskip
\noindent In all other cases, the the Esscher parameter does not exist. 
\e{Prop}
\proof In all cases (i) -- (iv), in view of Proposition \ref{prop:psi}, the monotonically increasing function $\psi$ is continuous on $[a,b]$ with range being an interval containing zero. Hence there exists $\kappa_0$ such that $\psi(\kappa_0)=0$. In all other cases, if $\kappa\in\mb R$ and $\psi(\kappa)$ exists, it can easily be seen that $\psi(\kappa)$ cannot be equal to zero.\hfill$\Box$
\bg{Rems}\label{Rems-EP}
{\rm (i)} If $E=\emptyset$, then $\psi(\kappa)$ does not exist or $\psi(\kappa)$ is equal to $-\infty$ or $+\infty$, for every $\kappa\in\mb R$. Hence there does not exist the Esscher parameter.  

\smallskip
{\rm (ii)} Condition (iv) of Proposition \ref{Esscher Parameter} includes the case that $a=b=0$. Then it follows $\psi(0)=0$, meaning that $\mb E\left[\xi\right]=0$, the Esscher parameter is $\kappa_0=0$, and the Esscher martingale transform is $\mb P$. 
\e{Rems}

}
\end{appendix}

\begin{thebibliography}{99}
{\small
\bibitem{AA} A. Andrusiv:
\textit{Stochastic Models for Finance and Insurance. Relative
Entro\-py.} PhD Thesis, University of Jena 2017, available online: \begin{small}
https://www.db-thueringen.de/servlets/MCRFile
\end{small}
\begin{small}
NodeServlet/dbt\_derivate\_00038814/Thesis.pdf
\end{small}
\bibitem{AE} A. Andrusiv, H.-J. Engelbert: On the minimal entropy martingale measure for L\'evy processes. To appear in: \textit{Stochastics. An International Journal of Probability and Stochastic Processes.} 
\bibitem{Ca-Th10} S. Cawston:
\emph{Mod\`{e}les de L\'evy exponentiels en finance:
mesures de f-divergence minimale et mod\`{e}les avec change-point.}
These de Doctorat. Math\'ematiques Appliqu\'ees, \'Ecole Doctoral STIM, Universit\'e d'Angers, 2010.
\bibitem{C-V-3} S. Cawston, L. Vostrikova: 
On continuity properties for option prices in exponential L\'evy models.
\newblock \textit{Theory of Probability and Their Applications,} 54, no. 4, 645--670 (2009).
\bibitem{C-V-1} S. Cawston, L. Vostrikova:
L\'evy preservation and associated properties for f-divergence minimal equivalent martingale measures. In: Shiryaev A., Presman E., Yor M.: \textit{Prokhorov and Contemporary Probability Theory}, Springer-Verlag, pp. 163--197 (2013).
\bibitem{C-V-2} S. Cawston, L. Vostrikova: 
An f-divergence approach for optimal portfolios in exponential L\'evy models. In: Kabanov Yu., Zariphopoulou T., Rutkowski M.: \textit{Inspired by Finance}, Springer-Verlag, pp. 83--101 (2014).
\bibitem{C99} T.\ Chan: \ Pricing contingent claims on stocks driven by L\'evy processes. Annals of Applied Probability 9, 504--528 (1999).
\bibitem{CS02} A.\ Cherny, A. N.\ Shiryaev: On minimization and maximization of entropy in various disciplines. \emph{Theory Probab. Appl.,} \textbf{48}(3), 447--464 (2003).
\bibitem{Cs75} I.\ Csisz\'ar: \ I-divergence geometry of probability distributions and minimization problems. \emph{Ann. Prob.} \textbf{3}, 146--158 (1975).
\bibitem{CT03} R.\ Cont, P.\ Tankow: \emph{Financial modelling with Jump Processes}. Chapman \& Hall CRC Press, 2003.
\bibitem{Six-02} F.\ Delbaen, P.\ Grandits, T.\ Rheinl\"ander, D.\ Samperi, M.\ Schweizer, C.\ Stricker: \ Exponential hedging and entropic penalties. \emph{Mathematical Finance,} \textbf{12} (2), 99--123 (2002).
\bibitem{ES05} F.\ Esche, M.\ Schweizer:
Minimal entropy preserves the L\'evy property: \ how and why? 
\emph{Stochastic Process. Appl.,} \textbf{115}(2), 299--327 (2005).
\bibitem{E32} F.\ Esscher:
On the probability function in the collective theory of risk. \emph{Skand. Aktuarie Tidskr.,} \textbf{15}, 175--195 (1932).
\bibitem{Fr00} M.\ Frittelli: \ The minimal entropy martingale measure and the valuation problem in incomplete markets. \emph{Math.\ Finance,}\ \textbf{10}: 39--52 (2000).
\bibitem{FM03} T.\ Fujiwara, Y.\ Miyahara (2003): \ The minimal entropy martingale measures for geometric L\'evy processes.
\emph{Finance and Stochastics} \textbf{7}, 509--531 (2003).
\bibitem{GS94} H.U.\ Gerber, E.S.W.\ Shiu:
Option pricing by Esscher transforms. 
\emph{Trans.Soc.Actuar.,} \textbf{XLVI}, 98--140 (1994).
\bibitem{GS96} H.U.\ Gerber, E.S.W.\ Shiu:
Martingale approach to pricing perpetual American options
on two stocks. \emph{Mathematical Finance}, Vol \textbf{6}, No 3, 301--322 (1996).
\bibitem{HWY92} S.\ He, J.\ Wang, J.\ Yan: \emph{Semimartingale Theory and Stochastic Calculus.} Science Press, Beijing-New York, CRC Press INC. 1992.
\bibitem{HS06} F.\ Hubalek, C.\ Sgarra:
Esscher transforms and the minimal entropy martingale measure for exponential L\'evy models. \emph{Quantitative Finance,} Vol. \textbf{6}, No. 2, 125--145 (2006).
\bibitem{IkWa}
N. Ikeda, S. Watanabe, {\em Stochastic Differential Equations and Diffusion Processes.}\/ Tokyo: North-Holland Publ. 1989.
\bibitem{JS00} J.\ Jacod, A.\ Shiryaev: \emph{Limit Theorems for Stochastic Processes 2ed}. Springer, 2000.
\bibitem{KS02} Y.\ Kabanov, C.\ Stricker: \ 
On the optimal portfolio for the exponential utility maximization: remarks to the six-author paper. 
\emph{Mathematical Finance,} \textbf{12} (2), 125--134 (2002).
\bibitem{Sa04} K.\ Sato: \emph{L\'evy processes and infinitely divisible distributions}. Cambridge University Press 1999.
}
\end{thebibliography}
\end{document}